\documentclass[12pt,oneside,english]{amsart}
\usepackage[latin1]{inputenc}
\usepackage{amssymb}

\makeatletter
\newtheorem{theorem}{Theorem}
\newtheorem{lemma}{Lemma}

\newtheorem{remark}{Remark}

\usepackage{babel}

\makeatother
\begin{document}
\baselineskip=17pt

\title[Strengthening of Newman phenomenon]
{On monotonic strengthening of Newman-like phenomenon on
$(2m+1)$-multiples in base $2m$}

\author{Vladimir Shevelev}
\address{Departments of Mathematics \\Ben-Gurion University of the
 Negev\\Beer-Sheva 84105, Israel. e-mail:shevelev@bgu.ac.il}

\subjclass{11A63.}

\begin{abstract}
We obtain exact and asymptotic expressions for the excess of
$(2m+1)$-multiples with even digit sums in the base $2m$ on interval
$[0,(2m)^k)$.
\end{abstract}

\maketitle
\section{Introduction}
    Consider a Newman-like digit sum

\begin{equation}\label{1}
S_{2m+1,0,2m}(x)=S^{(m)}(x)=\sum_{0\leq n < x:n\equiv
0mod(2m+1)}(-1)^{\sigma(n)}
\end{equation}

where $\sigma(n)=\sigma_{2m}(n)$ is the sum of digits of $n$ in the
representation of $n$ in the base $2m$.

We prove the following results

\begin{theorem}\label{t1}

\begin{equation}\label{2}
S^{(m)}((2m)^k)=\begin{cases}
\frac{2}{2m+1}\sum^m_{l=1}\left(\tan\frac{\pi
l}{2m+1}\right)^k,\;\;if \; k\; is\;
even\\\frac{2}{2m+1}\sum^m_{l=1}\left(\tan\frac{\pi
l}{2m+1}\right)^k \sin\frac{2\pi l}{2m+1},\;\; if \; k \; is \;
odd\end{cases}.
\end{equation}
\end{theorem}

\begin{theorem}\label{t2}

For an $m\in\mathbb{N},\;\;x_k=x_{m,k}=(2m)^k$, we have

\begin{equation}\label{3}
S^{(m)}(x_k)\sim \begin{cases} \frac{2}{2m+1}x^\alpha_k,\; if\; k\;
is\; even \\\frac{2}{2m+1}(\sin\frac{\pi}{2m+1})x^\alpha_k,\;\; if
\;k\; is\; odd\;(k\rightarrow\infty) \end{cases}
\end{equation}

where
\end{theorem}

\begin{equation}\label{4}
\alpha=\alpha_m=\frac{\ln\tan\frac{m\pi}{2m+1}}{\ln
(2m)}=\frac{\ln\cot(\frac{\pi}{4m+2})}{\ln(2m)}.
\end{equation}

In particular, in the case of $m=1$ we obtain the Gelfond-Newman
constant $\alpha_1=\frac{\ln 3}{\ln 4}$ (see \cite{1}, \cite{3}).
Further, we have $\alpha_2= 0.81 0922\ldots,\; \alpha_3=
0.824520\ldots,\;\alpha_4= 0.834558\ldots, \alpha_5=
0.842306\ldots$. Note, that

\newpage

\begin{equation}\label{5}
1-\frac{\ln\frac \pi 2}{\ln(2m)}< \alpha_m<1
\end{equation}

and $\alpha_m$ tends to 1. Indeed, we have

$$
\sin\frac{\pi}{8m+4}>\frac{\sqrt{2}}{2}\cdot\frac 4
\pi\cdot\frac{\pi}{8m+4}=\frac{\sqrt{2}}{4m+2}
$$
and
$$
\sin\frac{\pi}{2m+2}>\frac{\sqrt{2}}{2m+1}.
$$

Therefore,

$$
\frac{4m}{\pi}<\frac{4m+2}{\pi}\left(1-\frac{1}{(2m+1)^2}\right)
<\cot\frac{\pi}{4m+2}<\sqrt{2}(m+1)
$$

and (\ref{5}) follows from (\ref{4}).

Theorem \ref{t2} allows as in \cite{4} to obtain the sharp
estimations for $S^{(m)}(x)$.

\section{Lemmas}

\begin{lemma}\label{l1}
$$
S^{(m)}((2m)^k)=\frac{1}{2m+1}\sum^{2m}_{l=0}\prod^{k-1}_{j=0}
\left(1-\omega^{l(2m)^j}+\omega^{2l(2m)^j}-\ldots
-\omega^{(2m-1)l(2m)^j}\right),
$$
where
$$
\omega=\omega_{2m+1}= e^{\frac{2\pi i}{2m+1}}.
$$
\end {lemma}

\slshape Proof. \upshape The right hand side evidently equals to

$$
\frac{1}{2m+1}\sum^{2m}_{l=0}\sum^{(2m)^k-1}_{r=0}(-1)^{\sigma(r)}\omega^{lr}=
\frac{1}{2m+1}\sum^{(2m)^k-1}_{r=0}(-1)^{\sigma(r)}\sum^{2m}_{l=0}(\omega^l)^r=
$$
$$
=\frac{1}{2m+1}\sum^{(2m)^k-1}_{r=0}(-1)^{\sigma(r)}\cdot\begin{cases}0,\;
if\; r\; is \; not \; a\; multiple \; of \; 2m+1\\
2m+1,\; if \; r \; is \; a \; multiple \; of \; 2m+1\end{cases}=
$$
$$
\sum^{(2m)^k-1}_{r=0,(2m+1)|
r}(-1)^{\sigma(r)}=S^{(m)}((2m)^k).\blacksquare
$$

\begin{lemma}\label{l2}
$$
\frac{\omega^l-1}{\omega^l+1}=i\tan\frac{\pi
l}{2m+1},\;\;l=1,2,\ldots,2m.
$$
\end{lemma}

\slshape Proof. \upshape  Straightforward, after some quite
elementary transformations.$\blacksquare$

\newpage

\section{Proof of Theorem 1}

Using Lemma 1 and noticing that

$$
(2m)^j\equiv (-1)^j(mod(2m+1))
$$

we have

$$
S^{(m)}((2m)^k)=
$$
$$
=\begin{cases} \frac{1}{2m+1}\sum\limits^{2m}_{l=0}\left(
\left(\sum\limits^{2m-1}_{j=0}(-1)^j\omega^{lj}\right)^{\frac k
2}\left(\sum\limits^{2m-1}_{j=0}(-1)^j\omega^{-lj}\right)^\frac k
2\right),\;if\; k\; is \; even\\
\frac{1}{2m+1}\sum\limits^{2m}_{l=0}\left(\left(\sum\limits^{2m-1}_{j=
0}(-1)^j\omega^{lj}\right)^{\frac{k+1}{2}}
\left(\sum\limits^{2m-1}_{j=0}(-1)^j\omega^{-lj}\right)^\frac{k-1}{2}\right),\;
if \; k \; is \; odd\end{cases}=
$$
$$
=\begin{cases}\frac{1}{2m+1}\sum\limits^{2m}_{l=0}\left(\frac{1-\frac{1}{\omega^l}}
{1+\omega^l}
\cdot\frac{1-\omega^l}{1+\frac{1}{\omega^l}}\right)^{\frac
k 2},\;if\; k\; is\; even\\
\frac{1}{2m+1}\sum\limits^{2m}_{l=0}\left(\frac{1-\frac{1}{\omega^l}}{1+\omega^l}
\cdot\frac{1-\omega^l}{1+\frac{1}{\omega^l}}\right)^\frac{k-1}{2}
\left(\frac{1-\frac{1}{\omega^l}}{1+\omega^l}\right), \;if\; k\;
is\; odd\end{cases}=
$$
$$
=\begin{cases}\frac{1}{2m+1}\sum\limits^{2m}_{l=1}(-1)^\frac k
2\left(\frac{1-\omega^l}{1+\omega^l}\right)^k,\; if\; k \; is\;
even\\\frac{1}{2m+1}\sum\limits^{2m}_{l=1}(-1)^\frac{k+1}{2}
\left(\frac{1-\omega^l}{1+\omega^l}\right)^k\omega^{-l},\;
if\;k\;is\;odd\end{cases}.
$$

Using Lemma 2 we find

$$
S^{(m)}\left((2m)^k\right)=\begin{cases}\frac{1}{2m+1}\sum\limits^{2m}_{l=1}
\left(\tan\frac{\pi l}{2m+1}\right)^k,\;if\; k \;is \; even\\
\frac{1}{2m+1}\sum\limits^{2m}_{l=1} \left(\tan\frac{\pi
l}{2m+1}\right)^k i\omega^{-l},\;if\; k \;is \; odd.\end{cases}
$$

Finally, notice that for $l=1,2,\ldots,m$ we have

$$
\tan\frac{\pi l}{2m+1}=-\tan\frac{\pi(2m-(l-1))}{2m+1},
$$
$$
\left(\tan\frac{\pi
l}{2m+1}\right)i\omega^{-l}-\tan\frac{\pi(2m-(l-1))}{2m+1}i\omega^{-(2m-(l-1))}=
$$

$$
=i\left(\tan\frac{\pi
l}{2m+1}\right)\left(\omega^{-l}-\omega^l\right)=2\tan\frac{\pi
l}{2m+1}\cdot\sin\frac{2\pi l}{2m+1},
$$

and the theorem follows $\blacksquare$

\newpage

\section{Proof of Theorem 2}

Choosing the maximal exponent in (\ref{2}) we find for
$k\rightarrow\infty$

$$
S^{(m)}\left((2m)^k\right)\sim \begin{cases}
\frac{2}{2m+1}\left(\tan\frac{m\pi}{2m+1}\right)^k,\; if\; k\; is\;
even\\
\frac{2}{2m+1}\sin\frac{2\pi
m}{2m+1}\left(\tan\frac{m\pi}{2m+1}\right)^k,\; if\; k\; is\;
odd\end{cases}
$$

It is left to notice that if

$$
\ln S^{(m)}\left((2m)^k\right)\sim
k\ln\tan\frac{m\pi}{2m+1}=\ln\left((2m)^{k\alpha}\right)
$$

then the number $\alpha=\alpha_m$ is defined by (\ref{4}).
$\blacksquare$

\begin{remark} 1. Notice that, $S^{(m)}(2m)=1$. For $k\geq 2$, using
an equivalent representation
\end{remark}
$$
S^{(m)}\left((2m)^k\right)=\begin{cases}\frac{2}{2m+1}\sum\limits^{m-1}_{\lambda=0}
\left(\cot\frac{(2\lambda+1)\pi}{4m+2}\right)^k,\;if\; k\; is\;
even\\
\frac{2}{2m+1}\sum\limits^{m-1}_{\lambda=0}
\left(\cot\frac{(2\lambda+1)\pi}{4m+2}\right)^k\sin\frac{2\pi(m-\lambda)}{2m+1},\;
if\; k\; is\; odd\end{cases}
$$

and considering
$\sum\limits^{m^\sigma}_{\lambda=0}+\sum\limits^{m-1}_{\lambda=m^\sigma+1}$
with $0<\sigma< 1$, it is easy to obtain an interesting formula for
$m\rightarrow\infty$ (and a fixed $k\geq 2$)

$$
S^{(m)}\left((2m)^k\right)\sim\begin{cases}\frac{2^k(2^k-1)}{\pi^k}\zeta(k)m^{k-1},\;
if\;k\geq 2 \; is\; even,\\
\frac{2^k(2^{k-1}-1)}{\pi^{k-1}}\zeta(k-1)m^{k-2}, \; if\; k\geq 3\;
is\; odd,\end{cases}
$$

or in terms of the Bernoulli numbers,(see,e.g.\cite{6})

$$
S^{(m)}\left((2m)^k\right)\sim\begin{cases}2^{2k-1}(2^k-1)\frac{|B_k|}{k!}m^{k-1},
\; if\; k\geq 2\; is\; even,\\
2^{2k-2}(2^{k-1}-1)\frac{|B_{k-1}|}{(k-1)!}m^{k-2}, \; if\; k\geq
3\; is\; odd.\end{cases}
$$

E.g., for $k=4,\;\;|B_4|=\frac{1}{30}$, therefore,

$$
S^{(m)}\left((2m)^4\right)\sim \frac{2^7\cdot 15}{24\cdot
30}m^3=\frac 8 3 m^3,
$$
while exactly we have a polynomial

$$
S^{(m)}\left((2m)^4\right)=\frac{2m}{3}(4m^2+6m-1).
$$

\newpage

\begin{remark} If $k$ is even then $\frac{2m+1}{2}S^{(m)}((2m)^k)$
is the k-th power sum of the roots of the polynomial
\end{remark}
$$
\sum^m_{r=0}(-1)^r\begin{binom} {2m+1}{2r}\end{binom}x^{2m-2r}=0.
$$

Indeed, all roots of this polynomial are:  $\tan{\frac{\pi
l}{2m+1}},\;\;l=1,2,\ldots,m $. (cf.\cite{5}). Using the Littlewood
formula for the power sum in a determinant form \cite{2}, it is easy
to see that $\frac 1 2 S^{(m)}((2m)^k)$ is a polynomial of degree
$k-1$ with the integer values for $m\in\mathbb{N}$. Thus,
$S^{(m)}((2m)^k)$ always is even for even $k$.

     What could say about the case of an odd $k$?

\begin{remark}  As an additional corollary, for even $k$ we have
\end{remark}
$$
B_k=\frac{(-1)^{\frac k 2 -1}k!}{2^{2k}(2^k-1)}\det A,
$$

where $A=\{a_{ij}\}$ is $k\times k$ matrix with
$\frac{3k^2-4k+4}{4}$ zeros with the following nonzero elements:
 $a_{j,j+1}=1,\;j=1,2,\ldots,k-1;\;
a_{t,t-(2l-1)}=\frac{(-4)^l}{(2l-\delta_{l,\left\lfloor\frac{t+1}{2}\right\rfloor})!},
\;\;l=1,2,\ldots,\lfloor\frac t 2\rfloor,\; t=2,3,\ldots, k$, where
$\delta_{i,j}$ is the Kronecker symbol.

For example,

$$
B_4=-\frac{24}{256\cdot 15}\begin{vmatrix}0&  1&  0 & 0 \\-4 & 0& 1
&0\\ 0 &-2 & 0&  1 \\\frac 8 3& 0&-2& 0 \end{vmatrix}=\frac{1}{160}
\begin{vmatrix}-4&  1&  0 \\0& 0& 1
\\ \frac 8 3& -2& 0 \end{vmatrix}=
$$

$$
=-\frac{1}{160}\begin{vmatrix}-4 &  1 \\\frac 8 3&
-2\end{vmatrix}=-\frac{1}{160}\cdot\frac{16}{3}=-\frac{1}{30}.
$$

\;\;

\end{document}